\documentclass[11pt]{article}

\usepackage[utf8]{inputenc}

\usepackage{times,amssymb,amsmath,exscale,array,latexsym}
\usepackage{graphicx}
\graphicspath{{figs/}}         
\usepackage{epsfig}
\usepackage{color}
\usepackage[dvipsnames]{xcolor}
\usepackage{subcaption}
\definecolor{marin}{rgb}   {0.,   0.3,   0.7}
\definecolor{rouge}{rgb}   {0.8,   0.,   0.}
\definecolor{sepia}{rgb}   {0.8,   0.5,   0.}
\usepackage[colorlinks,citecolor=marin,linkcolor=rouge,
            bookmarksopen,
            bookmarksnumbered
           ]{hyperref}

\newcommand{\R}{\mathbb{R}}

\newcommand{\e}{\ensuremath{\mathrm{e}}}

\numberwithin{equation}{section}

\newcommand{\QED}{\mbox{}\hfill \raisebox{-0.2pt}{\rule{5.6pt}{6pt}\rule{0pt}{0pt}}
          \medskip\par}

\newcommand{\cL}{{\mathcal L}}

\makeindex
\begin{document}

\title{Generalized extrapolation methods based on compositions of a basic 2nd-order scheme}

\author{S. Blanes$^{1}$, F. Casas$^{2}$, L. Shaw$^{3}$ \\[2ex]
$^{1}$ {\small\it Universitat Polit\`ecnica de Val\`encia, Instituto de Matem\'atica Multidisciplinar, 46022-Valencia, Spain}\\{
\small\it email: serblaza@imm.upv.es}\\[1ex]
$^{2}$ {\small\it Departament de Matem\`atiques and IMAC, Universitat Jaume I, 12071-Castell\'on, Spain}\\{
\small\it email: Fernando.Casas@mat.uji.es}\\[1ex]
$^{3}$ {\small\it Departament de Matem\`atiques, Universitat Jaume I, 12071-Castell\'on, Spain}\\{
\small\it email: shaw@uji.es}\\[1ex]
}


%
\maketitle

\begin{abstract}

We propose new linear combinations of compositions of a basic second-order scheme with appropriately chosen coefficients to construct higher order
numerical integrators for differential equations. They
can be considered as a generalization of extrapolation methods and multi-product expansions. A general analysis is provided and new methods up to order
8 are built and tested. The new approach is shown to reduce the latency problem when implemented in a parallel environment and leads to schemes that are significantly more
efficient than standard extrapolation when the linear combination is delayed by a number of steps.

\end{abstract}\bigskip



\section{Introduction}
\label{sect1}

Extrapolation methods constitute a class of efficient numerical integrators for the initial value problem 
\begin{equation} \label{eq.1.ODE}
		x'=f(x),  \qquad  x(t_0)=x_0, \qquad\qquad x \in \mathbb{R}^d,
\end{equation}
especially when high accuracy is desired and relatively short time integrations are considered. Although it is possible to construct extrapolation methods
starting from a basic first-order scheme, better results are achieved in general with a time-symmetric 2nd-order scheme $S_{h}$,
since the asymptotic expansion of its local error contains only even powers of the step size $h$ \cite{hairer93sod}. If $\varphi_t$
denotes the flow of (\ref{eq.1.ODE}), i.e., $x(t) = \varphi_t(x_0)$, then
\[
  x_{n+1} = S_h(x_n) = \varphi_h(x_n) + \mathcal{O}(h^3),
\]  
where $x_n$ is taken as the approximation of $x(t_n)$ at $t_n = t_0 + n h$. 
 Given a sequence of integer numbers $0<m_1<m_2< m_3 < \ldots$, then it is possible to achieve order $2r$, $r=2,3,\ldots$, by determining the
 coefficients $b_i$ in the \emph{multi-product expansion} (MPE) \cite{chin11mpo}
\begin{equation}\label{eq.extrap} 
  \psi_h^{[2r]} = \sum_{i=1}^r b_i  \, \big( S_{h/m_i} \big)^{m_i}, \qquad\quad r = 2, 3, \ldots
\end{equation}
 so that
\[
  \psi_h^{[2r]}(x_0) = \varphi_h(x_0) + \mathcal{O}(h^{2r+1)}.
\]
Here $\big( S_{h/2} \big)^{2} \equiv S_{h/2}^{}\circ S_{h/2}^{}$, etc.
In particular, if one takes the harmonic sequence $m_i = i$, the resulting linear combination
\begin{equation}\label{h_s1} 
  \psi_h^{[2r]} = \sum_{i=1}^r b_i  \, \big( S_{h/i} \big)^{i}
\end{equation}
produces an extrapolation method of order $2r$ with only $r(r+1)/2$ evaluations of the basic map $S_h$
\cite{hairer93sod,chin11mpo}. For orders 4, 6 and 8 one has, respectively
\begin{equation} \label{em46}
\begin{aligned}
  & \psi_h^{[4]} = -\frac13 S_{h} + \frac43 S_{h/2}^{2}  \\
  & \psi_h^{[6]} = \frac{1}{24} S_h - \frac{16}{15} S_{h/2}^2 + \frac{81}{40} S_{h/3}^3 \\
  & \psi_h^{[8]} = -\frac{1}{360} S_h + \frac{16}{45}  S_{h/2}^2 - \frac{729}{280} S_{h/3}^3 + \frac{1024}{315} S_{h/4}^4.
\end{aligned}
\end{equation}
In fact, analytic expressions for the coefficients $b_i$ exist at any order \cite{chin10mps}, so that construction of MPEs / extrapolation methods of arbitrarily high order is trivial.

Schemes (\ref{h_s1}) are also well suited for parallelization: each processor can evaluate the composition $\big( S_{h/i} \big)^{i}$, requiring at most
$r$ evaluations of the basic map $S_h$. There remains, however, the problem of latency: the period of time used for
communication between processes during which none of the
processes may advance. In the application under consideration, latency will take place whilst one sums the outputs of the various compositions. This is
particularly problematic for schemes (\ref{eq.extrap}), since all processes require different numbers of evaluations and so all have to wait until the composition involving
the largest number of maps terminates before proceeding further. As a result, communication between processors may be more expensive
than computing the different compositions in parallel.

It is then of interest to examine the case where one delays the summation by, say, $p > 1$ steps, i.e., when in place of (\ref{eq.extrap}) one uses
\begin{equation}\label{eq.extrap_p}
  \hat{\psi}_{ph}^{[2r]} = \sum_{i=1}^r b_i \left( S_{h/m_i}\right)^{m_i p}
\end{equation} 
to reduce the communication between processes. One should stress, however, that $\hat{\psi}_{ph}^{[2r]} \ne \left(\psi_h^{[2r]}\right)^p$, so that additional
errors depending on $p$ and even order reductions may take place when one uses very large values of $p$ in (\ref{eq.extrap_p}).
Notice, however, that even in this case not all processes involve the same amount of computations, and so most of  them have to wait until the processor
evaluating the largest composition finishes.

Multi-product expansions of the form (\ref{eq.extrap}) have shown to be competitive with other integration schemes on a number of problems, especially
when high accuracy is required \cite{chin10mps,chin11mpo}, even when a fixed step size is used. On the other hand, they do not belong to the class
of structure-preserving methods. Specifically, suppose problem (\ref{eq.1.ODE}) is formulated in some Lie group so that the exact solution preserves some qualitative features related with this fact. For instance, if (\ref{eq.1.ODE}) corresponds to a Hamiltonian system, then the exact flow is symplectic. In that
case, the numerical solution furnished by an MPE of order $2r$ is no longer symplectic, but only up to one order higher than the order of the method itself. If (\ref{eq.1.ODE}) 
evolves in the SU($N$) group, then the exact matrix solution is unitary (with unit determinant), whereas the numerical approximation preserves unitarity up up to
order $2r+1$, etc.
Moreover,
MPEs may be more prone to round-off errors than other
types of integrators, due to the different size of the terms to be added up.

Having identified some of the advantages and disadvantages of MPEs, it is natural to generalize them so as to remove their drawbacks whilst 
still maintaining their characteristic good efficiency. A step in that direction is
provided by the study carried out in \cite{blanes21npi}, 
where new 4th- and 6th-order schemes show greater efficiency than MPEs
 (\ref{em46}) on the 2-dimensional Kepler problem. The schemes proposed in  \cite{blanes21npi} can be formulated as particular instances of the general linear combination
\begin{equation}\label{linear_comb1}
  \psi_h^{[k,m]} = \sum_{i=1}^k b_i \prod_{j=1}^{m} S_{a_{ij}h},
\end{equation}
where, as before, $S_h$ is a 2nd-order time-symmetric integrator and the condition $\sum_{i=1}^k b_i = \sum_{j=1}^m a_{ij} = 1$ 
is usually taken for consistency.

The presence of additional coefficients in (\ref{linear_comb1}) may be used for various purposes. In particular, one could increase the order of
preservation of whatever properties system (\ref{eq.1.ODE}) may possess (such as symplecticity or unitarity), or to reduce the most significant contributions to the truncation error. In contrast,
MPEs use all available coefficients to increase the order of the method as much as possible, and in doing so, they provide excellent results, as we will
see in section \ref{sect2}. It turns out that there exist order barriers depending on the value of $m$, so that, even for arbitrarily large $k$ in 
\eqref{linear_comb1}, it is not possible to raise the order of the method. Nonetheless, one can still use this extra freedom to improve the efficiency (although not the order).

Each of the integrations in the sum \eqref{linear_comb1}
may be performed simultaneously and then combined, but in this case all the processes require the same amount of
work per step, so that in principle the latency period is reduced. If we are interested in further reductions, then we can combine the solutions after $p$ steps, i.e.,
in place of (\ref{linear_comb1}) we use instead
\begin{equation}\label{linear_comb_p}
  \psi_h^{[k,m,p]} = \sum_{i=1}^k b_i \left( \prod_{j=1}^{m} S_{a_{ij}h} \right)^p, \qquad\quad \sum_{i=1}^k b_i = \sum_{j=1}^m a_{ij} = 1.
\end{equation}
Here we carry out a general analysis of schemes (\ref{linear_comb1}) and (\ref{linear_comb_p}), eventually proposing new schemes that perform better
than extrapolation methods for a number of problems, whilst possessing better preservation properties and reducing the latency problem.

The structure of the paper is as follows. Section \ref{sect2} collects information previously available about multi-product expansions, 
in particular the order conditions and their general solutions, as well as the structure of the truncation error and their preservation properties. The general
analysis of linear combinations (\ref{linear_comb1}) is detailed in section \ref{sect3}, where we also construct new schemes within this family of orders 4, 6
and 8. The new methods are tested in {section \ref{sect4} }on two examples in comparison with multi-product expansions and the integrators proposed in
\cite{blanes21npi}. Finally, section \ref{sect5} analyzes the latency problem for this type of schemes and shows how methods preserving qualitative properties up to
higher orders can be successfully used when the summation of the different compositions is delayed by a number of steps. 

\section{Multi-product expansions}
\label{sect2}

\paragraph{Series of differential operators.}

The Lie formalism constitutes an appropriate tool for the analysis of the methods considered in this work. As is well known, associated with $f$ in the
ODE (\ref{eq.1.ODE}) there exists an operator $F$, called the Lie derivative, and defined by
\[
    F \, g (x) = \left.\frac{d}{d h}\right|_{h=0} g(\varphi_{h}(x)) 
\]
for each smooth function $g:\mathbb{R}^d \rightarrow \mathbb{R}$ and $x \in \mathbb{R}^d$, so that
\begin{equation}   \label{eq:3b}
(F\, g) (x) = f(x) \cdot \nabla g(x).
\end{equation}
Then, the $h$-flow of (\ref{eq.1.ODE}) {satisfies} \cite{hairer06gni,sanz-serna94nhp}
\[
   g(\varphi_{h}(x)) =\big( \e^{h F}g \big) (x) = \sum_{k=0}^{\infty} \frac{h^k}{k!} (F^k g)(x).
\]   
Analogously, one can associate a series of linear operators to the basic 2nd-order  method $S_h$ as
 \cite{blanes08sac}
\[
    g(S_h(x)) = \big(\e^{Y(h)}g \big)(x), 
\]
for all functions $g$, where
\[
  Y(h) = \sum_{k \ge 0} h^{2k+1} Y_{2k+1} \qquad \mbox{ and }  \qquad Y_1 = F.
\]  
 Notice that only odd powers of $h$ appear in $Y(h)$, due to the time-symmetry of $S_h$.

In that case, the composition $\big( S_{h/m_i} \big)^{m_i}$ in (\ref{eq.extrap}) has associated a series $S^{(m_i)}(h)$
of differential operators given by
\begin{equation} \label{series}
  S^{(m_i)}(h) = \e^{m_i Y(h/m_i)}  = \exp \left( \sum_{k = 0}^{\infty} \frac{h^{2 k + 1}}{m_i^{2 k}} \, Y_{2 k + 1} \right),
\end{equation}
whereas the series associated to the MPE (\ref{eq.extrap}) is the linear combination
\begin{equation} \label{lcop_1}
  \Psi(h) = \sum_{i = 1}^r  b_i \, S^{(m_i)}(h).
\end{equation}

\paragraph{Order conditions for MPEs.}

We then proceed as in \cite{blanes99eos} and write (\ref{lcop_1}) as
\[
   \Psi(h) = \e^{\frac{h}{2} Y_1} \, Z(h) \, \e^{\frac{h}{2} Y_1},
\]
where the expression of $Z(h)$ can be obtained with the Baker--Campbell--Hausdorff formula \cite{varadarajan84lgl} as
\begin{equation} \label{eq_Z}
  Z(h) =  \sum_{i = 1}^r b_i \, \e^{h^3 W_i(h)}, 
\end{equation}
with
\[
  W_i(h) = \frac{1}{m_i^2} \, R^{(1)}\left(\frac{h}{m_i}\right) = \frac{1}{m_i^2} \, 
  \sum_{j=0}^{\infty} \left(\frac{h}{m_i}\right)^{2j} R_{2j}^{(1)}, \qquad \mbox{ and } \qquad 
  R_{2j}^{(1)} = Y_{3+2j} + \sum_{\ell = 0}^{j-1} m_i^{2j-2 \ell} u_{\ell}.
\]  
Here $u_{\ell}$ is an element of the free Lie algebra $\mathcal{L}(Y_j)$
generated by the operators $\{ Y_{2 \ell +1} \}$, $\ell = 0,1,2, \ldots$. This can be considered as the vector space spanned by
all $Y_{2\ell +1}$ and their independent nested commutators \cite{munthe-kaas99cia}.

Expanding the
exponentials in (\ref{eq_Z}) we get
\begin{equation} \label{eq_Z2}
  Z(h) = G_0 I + \sum_{\ell = 1}^{\infty} \frac{1}{\ell!} \, h^{3 \ell} \, \sum_{j=0}^{\infty} h^{2j} \, G_{2 \ell + 2 j} \, R_{2j}^{(\ell)}
\end{equation}  
in terms of 
\[
  G_s = \sum_{i=1}^r \frac{b_i}{m_i^s}, \qquad\qquad R_{2s}^{(\ell)} = \sum_{i=0}^s R_{2i}^{(1)} R_{2s-2i}^{(\ell - 1)}, \qquad  \ell > 1, \quad s = 0, 1, 2, \ldots
\]
In more detail, 
\[
\begin{aligned}
 Z(h)  = & \, G_0 I + h^3 \left( G_2 R_0^{(1)} + h^2 G_4 R_2^{(1)} + h^4 G_6 R_4^{(1)} + h^6 G_8 R_6^{(1)} + \mathcal{O}(h^8) \right) \\
      & + \frac{1}{2} h^6 \left(G_4 R_0^{(2)} + h^2 G_6 R_2^{(2)} +\mathcal{O}(h^4) \right) \\
      & +  \frac{1}{6} h^9 \left(G_6 R_0^{(3)} + \mathcal{O}(h^2) \right) + \cdots
\end{aligned}
\]      
It is then clear that if
\begin{equation} \label{orcon1}
   G_0 = 1, \qquad\qquad G_{2 \ell} = 0, \quad \ell= 1, \ldots, r - 1,
\end{equation}
then one has a method of order $2r$, $r \ge 2$. The order conditions (\ref{orcon1}) form a linear system in the coefficients $b_i$, with unique solution
given by \cite{chin10mps}
\[
  b_i = \prod_{j = 1 \atop j \ne i}^r \frac{m_j^{-2}}{m_j^{-2} - m_i^{-2}} = \prod_{j = 1 \atop j \ne i}^r \frac{m_i^2}{m_i^2 - m_j^2}.
\]  
In fact, the structure of $Z(h)$ in (\ref{eq_Z2}) allows one to get additional information on the extrapolation method (\ref{eq.extrap}). Thus, in particular,
it is possible to compute analytically the main term in the truncation error. If $\Psi(h)$ in (\ref{lcop_1}) corresponds to a method of order $2r$, then the
term in $h^{2r+1}$ clearly comes from the term $\ell = 1$ in the series (\ref{eq_Z2})\footnote{The first non-vanishing term coming from $\ell = 2$ is proportional to
$h^{2r+2}$ and contributions from $\ell \ge 3$ contain higher powers of $h$.}, which reads
\[
  G_{2r} \, R_{2r-2}^{(1)} = G_{2r} \, Y_{2r+1} + \sum_{\ell = 0}^{r-2} G_{2+2\ell} \, u_{\ell} = G_{2r} \, Y_{2r+1}, 
\]
since $G_2 = \cdots =G_{2r-2} = 0$ for a method of order $2r$. In consequence, 
\[
   \e^F = \Psi(h) + h^{2r+1} G_{2r} Y_{2r+1} + \mathcal{O}(h^{2r+2}), 
\]
and $G_{2r}$ can be computed  as follows:
\begin{equation} \label{g2r}
\begin{aligned}
 G_{2r} &= \sum_{i=1}^r \frac{b_i}{m_i^{2r}} = \sum_{i=1}^r m_i^{-2r} \prod_{j = 1 \atop j \ne i}^r \frac{m_i^2}{m_i^2 - m_j^2}=\left(\prod_{j=1}^rm_j^{-2}\right)\sum_{i=1}^r\prod_{j=1\atop j \ne i}^r m_j^{2}\prod_{j=1\atop j \ne i}^r \frac{1}{m_i^2-m_j^2}\\
&=(-1)^{r-1}\left(\prod_{j=1}^rm_j^{-2}\right)\sum_{i=1}^r\prod_{j=1\atop j \ne i}^r \frac{m_j^2}{m_j^2-m_i^2}=
(-1)^{r-1} \, \prod_{j=1}^r \frac{1}{m_j^{2}},
\end{aligned}
\end{equation}
where we have used the well-known property of the Lagrange polynomials $\sum_{i=1}^r \prod_{\substack{j=1\\j\neq i}}^r x_j/(x_j-x_i)=1$ with $x_j=m_j^2$ \cite{chin10mps}.
Then 
\[
  Z(h) = I + h^{2r+1} U + \frac{1}{2} h^{2r+2} G_{2r} R_{2r-4}^{(2)} + \cdots
\]
with
\[
  U = \sum_{j=0}^{\infty} h^{2j} \, G_{2r+2j} \, R_{2(r-1)+2j}^{(1)}.
\]
The important point is that $U$ contains only operators in $\mathcal{L}(Y_j)$, so that 
\[
 Z(h) = \e^{h^{2r+1} U} +  \frac{1}{2} h^{2r+2} \, G_{2r} \, R_{2r-4}^{(2)} + \mbox{ higher order terms. }   
\]
 Since $R_{2r-4}^{(2)}$ does not belong to $\mathcal{L}(Y_j)$, it is clear that the method $\psi_h^{[2r]}$ of order $2r$ only preserves 
 geometric properties related with $\mathcal{L}(Y_j)$ up to order $2r+1$. In particular, if (\ref{eq.1.ODE}) corresponds to a Hamiltonian system, then the
 numerical approximation furnished by a
scheme (\ref{eq.extrap}) of order $2r$ is symplectic up to order $2r+1$. 

It is illustrative to compare the main error term $G_{2r}$, and therefore the efficiency, for different sequences of integers $\{ m_i \}$, $i=1,2,\ldots$,
 commonly used in
extrapolation, and in particular for the Romberg sequence $m_i = 2^{i-1}$, the harmonic sequence $m_i = i$ and the Burlirsch sequence
\[
  m_i = 1, 2, 3, 4, 6, 8, 12, 16, 24, 32, \ldots
\]
If we define the efficiency of the extrapolation method \eqref{eq.extrap} of order $2r$  as
\[
  \mathcal{E}_f = n_s \, |G_{2r}|^{1/(2r)}
\]
where $n_s$ is the total number of evaluations of the basic map $S_h$ required by the scheme, then a straightforward computation shows that, for any order
$2r$,
\[
    \mathcal{E}_f \text{ harmonic} \le \mathcal{E}_f \text{ Bulirsch} \le \mathcal{E}_f \text{ Romberg},
\]    
in accordance with  \cite[p. 226]{hairer93sod}.

\section{Analysis of general linear combinations}
\label{sect3}

In view of the favorable features that MPEs possess, it is natural to consider now the more general family of integrators (\ref{linear_comb1})
\begin{equation} \label{eq.3.1o}
  \psi_h^{[k,m]} = \sum_{i=1}^k b_i \prod_{j=1}^{m} S_{a_{ij}h}, \qquad \mbox{ with } \qquad \sum_{i=1}^k b_i = \sum_{j=1}^m a_{ij} = 1.
\end{equation}
 Scheme
(\ref{eq.extrap}) is then recovered by taking $k=r$ terms and coefficients {$a_{ij} = 1/m_i$} for all $j$ in (\ref{eq.3.1o}). Whereas
this simple choice for the $a_{ij}$ leads to linear order conditions in the $b_i$, the general treatment is more involved. In any case, we
can proceed as in Section \ref{sect2} and write the operator $\Psi^{[k,m]}(h)$ associated to (\ref{eq.3.1o}) as 
\[
   \Psi_h^{[k,m]} = \e^{\frac{h}{2} Y_1} \, Z(h) \, \e^{\frac{h}{2} Y_1}, \qquad\qquad Z(h) =  \sum_{i = k}^r b_i \, \e^{h^3 W_i(h)},
\]
where now
\[
  W_i(h) = w_{31}^{(i)} \, Y_3 + h w_{41}^{(i)} \, [Y_1, Y_3] + h^2 (w_{51}^{(i)} \, Y_5 + w_{52}^{(i)} \, [Y_1,[Y_1,Y_3]]) + \mathcal{O}(h^3),  
\]
and $w_{mn}^{(i)}$ are polynomials in the coefficients $a_{ij}$. Expanding $Z(h)$ in powers of $h$ results in
\begin{equation} \label{expan_Z}
\begin{aligned}
  & Z(h) = G_{00} \, I + h^3 G_{31} E_{31} + h^4 G_{41} E_{41}  + h^5 \, \sum_{\ell =1}^2 G_{5 \ell} \, E_{5 \ell} \\
  & \quad + h^6 \left( \sum_{\ell =1}^2 G_{6 \ell} E_{6 \ell} + \frac{1}{2} \tilde{G}_{63} E_{31}^2 \right) + 
   h^7 \left( \sum_{\ell =1}^4 G_{7 \ell} E_{7 \ell} + \frac{1}{2} \tilde{G}_{75} \{E_{31}, E_{41} \} \right)  \\
   & \quad + h^8 \left( \sum_{\ell =1}^5 G_{8 \ell} E_{8 \ell} + \frac{1}{2} \left( \tilde{G}_{86} E_{41}^2 + \tilde{G}_{87} \{ E_{31}, E_{51} \}  +
      \tilde{G}_{88} \{ E_{31}, E_{52} \} \right) \right)+ \mathcal{O}(h^9),
\end{aligned}   
\end{equation}
where $E_{ij}$ denote the elements of the basis of the free Lie algebra $\mathcal{L}(Y_j)$ collected in Table \ref{table_1},  the symbol $\{ \cdot, \cdot \}$
denotes the anti-commutator, i.e.,  $\{E_{31}, E_{41} \}  = E_{31} E_{41} + E_{41} E_{31}$, and finally
\begin{equation} \label{coefficients}
\begin{aligned}
  & G_{00} = \sum_{i=1}^k b_i,  \qquad G_{mn}  = \sum_{i=1}^k b_i \, w_{mn}^{(i)}, \qquad 
  \tilde{G}_{63} = \sum_{i=1}^k b_i \, \big(w_{31}^{(i)} \big)^2, \qquad  \tilde{G}_{75} = \sum_{i=1}^k b_i \, w_{31}^{(i)}  w_{41}^{(i)}, \\
  & \tilde{G}_{86} =  \sum_{i=1}^k b_i \, \big(w_{41}^{(i)} \big)^2, \qquad \tilde{G}_{87} = \sum_{i=1}^k b_i \, w_{31}^{(i)}  w_{51}^{(i)}, \qquad 
    \tilde{G}_{88} = \sum_{i=1}^k b_i \, w_{31}^{(i)}  w_{52}^{(i)}.
\end{aligned}    
\end{equation}  
With expansion (\ref{expan_Z}) at hand, we can now analyze the situation order by order and eventually construct new schemes.

\begin{table}[t]
\small
\begin{center}
  \begin{tabular}{cllll}
     & \multicolumn{3}{c}{}  \\    \hline
 & $E_{31}=Y_3$ &   &    &  \\   \hline
 & $E_{41}=[Y_1,E_{31}]$   &  &    &  \\   \hline
 & $E_{51} = Y_5$  & $E_{52}=[Y_1, E_{41}]$  &      &  \\   \hline
 & $E_{61}= [Y_1,E_{51}]$   & $E_{62}=[Y_1,E_{52}]$ &   & \\ \hline
 & $E_{71}= Y_7$ & $E_{72}=[Y_1, E_{61}]$ & $E_{73}=[Y_1, E_{62}]$  & $E_{74} = [E_{31},E_{41}]$ \\ \hline
 & $E_{8\ell}=[Y_1,E_{7\ell}]$  & $\ell=1,\ldots,4$ &  $E_{85}=[E_{31},E_{51}]$ & \\ \hline
\end{tabular}
\caption{First terms in the the particular basis of $\cL(Y_1,Y_3,Y_5,Y_7)$ taken in this work. }\label{table_1}
\end{center}
\end{table}

\paragraph{Order 4.}
By taking 2-stage compositions, i.e., $m=2$ in (\ref{eq.3.1o}), 
\begin{equation} \label{psi_or4}
  \psi_h^{[k,2]} = \sum_{i=1}^k b_i \, S_{(1-a_i) h} \circ S_{a_i h},
\end{equation}  
one has enough parameters to construct a method of order 4. The aim here is to find particular solutions so that the resulting scheme is competitive
with the MPE $\psi_h^{[4]}$ of (\ref{em46}). Notice that in this case one must satisfy
\begin{equation} \label{or_4_2}
  G_{00} = 1, \qquad G_{31} = G_{41} = 0,
\end{equation}
where, explicitly,
\begin{equation} \label{doubleW}
\begin{aligned}
  w_{31}^{(i)} = (1-a_i)^3 + a_i^3,   & \qquad\qquad w_{41}^{(i)} = \frac{1}{2} a_i (2a_i-1)(1-a_i) \\
  w_{51}^{(i)} = (1-a_i)^5 + a_i^5,  &  \qquad\qquad w_{52}^{(i)} = \frac{1}{12} a_i (2 a_i - 1)^2 (a_i -1) + \frac{1}{24} w_{31}^{(i)}.
\end{aligned}
\end{equation}  
A 5th-order scheme would require, in addition to (\ref{or_4_2}), that $G_{51} = G_{52} = 0$. It turns out, however,  that the following constraint involving
the coefficients $w_{mn}^{(i)}$ in (\ref{doubleW}) holds:
\[
  \frac{5}{2} \, w_{31}^{(i)} - 4 \, w_{51}^{(i)} + 60 \, w_{52}^{(i)} = 1,
\]
so that conditions $G_{51}=0$ and $G_{52} = 0$ cannot be simultaneously satisfied and therefore a 5th-order barrier exists. 
In conclusion, with
$m=2$ only order 4 is possible with this approach, and the possibilities of improvement are quite limited. Two different options present themselves:
\begin{itemize}
 \item If  $G_{52} = 0$, then the only solution with $k=2$ corresponds to the standard extrapolation method $\psi_h^{[4]}$ of (\ref{em46}),
 leading to $G_{51} = -\frac{1}{4}$. 
 \item If $G_{51}=0$, there are no real solutions with $k=2$, whereas with $k=3$ one gets the solution proposed in \cite{blanes21npi}, with
 $G_{52} = 1/60$.
\end{itemize}
Of course, other alternatives are possible. Thus, with $k=2$ in (\ref{psi_or4}), we choose the solution providing the smallest values of $G_{51}$ and
$G_{52}$, whilst
keeping the size of the coefficients reasonably small. This results in a method with $G_{51} \approx -0.2089$, and $G_{52} \approx 0.0027$, denoted 
as $\psi_h^{[2,2]}$ in section \ref{sect4}.

Increasing the number of processors to $k=3$ introduces additional free parameters that can be used for various purposes. In particular,
\begin{itemize}
 \item We fix $G_{51}=0$, as in \cite{blanes21npi}, but explore solutions with the $b_i$ coefficients as small as possible. Thus, we construct the 
new method $\psi_h^{[3,2]}$, with $(\max_i (b_i) - \min_i(b_i)) = 4.59$. This should be compared with the value
 $(\max_i (b_i) - \min_i(b_i)) = 16.68$ corresponding to the equivalent scheme in  \cite{blanes21npi}.
 \item In addition to the order conditions (\ref{or_4_2}), we can also vanish $\tilde{G}_{63}$ and $\tilde{G}_{75}$ in the expression (\ref{expan_Z}) of $Z(h)$. 
 Notice that these terms do not belong to $\mathcal{L}(Y_j)$, so that, if the original problem corresponds to a Hamiltonian system, the resulting 4th-order method
 $\psi_{h,s}^{[3,2]}$ preserves symplecticity up to order 7 (instead of order 5, as for the other methods). In that case, we say that  $\psi_{h,s}^{[3,2]}$
 is \emph{pseudo-symplectic} of order 7 \cite{aubry98psr,casas21cop}.
\end{itemize} 

\paragraph{Order 6.}
From the previous discussion, it is clear that compositions with at least $m=3$ stages are required in (\ref{eq.3.1o}) to achieve order 6. We first
explore palindromic compositions, i.e., 
\begin{equation} \label{psi_or6}
  \psi_h^{[k,3]} = \sum_{i=1}^k b_i \, S_{a_i h} \circ S_{(1-2a_i) h} \circ S_{a_i h}.
\end{equation}  
From the expression of $Z(h)$ in (\ref{expan_Z}) it is clear that a 6th-order scheme only requires the following conditions to be satisfied:
\begin{equation} \label{oc_order6}
  G_{00} = 1, \qquad G_{31} = G_{51} = G_{52} = \tilde{G}_{63} = 0,
\end{equation}
since all terms $G_{ij}$ with even $i$ vanish by the symmetry of the composition. Moreover, from the expression of $w_{mn}^{(i)}$, not all terms are
independent. In particular, one has the following relations: 
\begin{equation} \label{constraints_6}
 \begin{aligned} 
  & G_{72} = -\frac{5}{4536} + \frac{5}{126} G_{71}, \qquad\quad G_{73} =  -\frac{1}{90720} + \frac{1}{2520} G_{71}, \\
  & G_{74} = -\frac{1}{648} + \frac{1}{18} G_{71}, \qquad\quad  
   -\frac{1}{15} \tilde{G}_{87} + 2 \tilde{G}_{88} + \frac{1}{30} G_{71} + G_{72} + 6 G_{73} - \frac{3}{2} G_{74} = 0,
 \end{aligned}
\end{equation} 
so that a 7th-order barrier exists.  If $G_{71}=0$, then $G_{72}$, $G_{73}$ and $G_{74}$ are fixed and moreover it is not possible to
simultaneously vanish the terms $\tilde{G}_{87}$ and $\tilde{G}_{88}$. On the other hand, if $G_{72} = G_{73} = G_{74} = 0$, then $G_{71} = 1/36$.  The extrapolation method $\psi_h^{[6]}$ of eq. (\ref{em46}) corresponds to this case when $k=3$, although other possibilities exist: in particular, we have
constructed a method which gives $G_{71} \approx 0.0199$ with small values of $b_i$, but it behaves in practice as $\psi_h^{[6]}$, so that it is necessary
to consider larger values of $k$. We then propose the following alternatives:
\begin{itemize}
 \item A $k=4$ method which vanishes $\tilde{G}_{87}$ and $\tilde{G}_{88}$, and thus it is pseudo-symplectic of order 8.
 \item Very often, the most important error term at order 7 corresponds to $G_{71}$. For this reason, we construct a $k=4$ method which vanishes $G_{71}$
 and $\tilde{G}_{87}$.
 \item The method presented in \cite{blanes21npi} with $k=5$ vanishes $G_{71}$, $\tilde{G}_{87}$ and $G_{91}$ (associated with $Y_9$ in the expansion of $Z(h)$).
 The dependences within the equations are such that these conditions fix all other terms at $h^7$, $h^8$ and $h^9$.
 \item A $k=5$ method which vanishes all the terms not belonging to $\mathcal{L}(Y_j)$ up to order 9, i.e., $\tilde{G}_{87}$, $\tilde{G}_{88}$ and
 $\tilde{G}_{99} \equiv  \sum_{i=1}^k b_i \, \big(w_{31}^{(i)} \big)^3$. We expect this additional preservation of symplecticity (for Hamiltonian systems) up to order
 9 to lead to better performance in the long run, although the fact that $G_{71} = 13/90$ may be a hindrance.
 \item Finally, a $k=5$ method vanishing $G_{71}$, $\tilde{G}_{87}$ and $G_{91}$ (as does the scheme in \cite{blanes21npi}), but with smaller values of the $b_i$ coefficients.
\end{itemize} 
We can also analyze  linear combinations $\psi_h^{[k,3]}$ involving non-palindromic compositions of the form 
\begin{equation} \label{npc6}
    S_{a_{i1}h} \circ S_{a_{i2}h} \circ S_{(1-a_{i1}-a_{i2})h},
\end{equation}   
 so that more free parameters are available for optimization. It turns out, however, that a
7th-order barrier is still present, since the following identity is satisfied: 
\[
  \frac{35}{4} G_{31} - \frac{63}{2} G_{51} + 35\tilde{G}_{63}+36G_{71}-420G_{72}+10080G_{73}-420G_{74}=1.
\]
It is thus possible to vanish all terms at order 7 except $G_{73} = 1/10080$, although the resulting method with $k=4$ is not particularly efficient in practice
(see Figure \ref{fig:6asymmetric}).

\paragraph{Order 8.}
We need compositions with at least $m=4$ stages in (\ref{eq.3.1o}). If we take symmetric compositions
\[
  S_{a_i h} \circ S_{(\frac{1}{2} - a_i) h} \circ S_{(\frac{1}{2} - a_i) h} \circ S_{a_i h},
\]
then 11 order conditions are required to achieve order 8, and thus $k=6$ terms provide enough parameters. It turns out, however, that any $\psi_h^{[k,4]}$
of this form is subject to a fifth-order barrier, since the identity
\[
  20 G_{31} - 144 G_{51} + 80 \tilde{G}_{63} = 1
\]
holds. If 5-stage symmetric compositions are considered, i.e., with 
linear combinations of the form
\begin{equation} \label{psi_or8}
  \psi_h^{[k,5]} = \sum_{i=1}^k b_i \, S_{a_{i1} h} \circ S_{a_{i2} h} \circ S_{(1-2a_{i1}-2 a_{i2}) h} \circ  S_{a_{i2} h} \circ S_{a_{i1} h},
\end{equation}  
with $k=4$ we are able to get order 8. One such method, which also eliminates the term $G_{91}$ is illustrated in
Section \ref{sect4} in comparison with  the extrapolation method $\psi_h^{[8]}$ (see Figure \ref{fig:5stage}).

On the other hand, with 4-stage non-symmetric compositions $ S_{a_{i1} h} \circ S_{a_{i2} h} \circ S_{a_{i3} h} \circ S_{a_{i4} h}$
the number of order conditions is 21. Although the extrapolation method $\psi_h^{[8]}$ in  (\ref{em46}) constitutes a particular solution with only $k=4$ terms,
in general we will need $k=7$ processors and the analysis is much more involved. 

The coefficients of all methods discussed here and tested in the succeeding section have been uploaded as \texttt{.txt} files to GitHub\footnote{\url{https://github.com/lshaw8317/Generalized-extrapolation-methods}}. In addition, the coefficients of the most efficient schemes are collected in the appendix.

\section{Numerical examples}
\label{sect4}
We next present some numerical experiments to illustrate the behavior of the new schemes described in Section \ref{sect3} on two 
simple problems of a different character, but both possessing conserved quantities: 
the 2-body gravitational problem and the Lotka--Volterra system. We will use the compensated sum technique
described in subsection \ref{sec:Rounding} to minimize roundoff error. It turns out that, with this technique and for the experiments we have carried out, methods with larger $b_i$ coefficients do not perform any differently to methods with smaller $b_i$, assuming they both satisfy the same order conditions.

\paragraph{Example 1: Kepler problem.} 
The motion of two bodies attracting each other through the gravitational law takes place in a plane and can be described by the 2-degrees-of-freedom
Hamiltonian 
\[
H(q,p)=T(p)+V(q)=\frac{1}{2}p^Tp- \mu \frac{1}{r},\qquad \mbox{ with }  \quad r = \|q\| = \sqrt{q_1^2 + q_2^2},
\]
$\mu = G M$, $G$ is the gravitational constant and $M$ is the sum of the masses of the two bodies.
We could take as the basic second-order integrator $S_h$ the St\"ormer--Verlet scheme \cite{hairer03gni}
\[
 S_h(q,p)=\phi_{h/2}^V\circ\phi_{h}^T\circ\phi_{h/2}^V(q,p), 
\]
 which corresponds to the symmetric composition  of the flows corresponding to the kinetic and potential energies with 
\[
\phi_{h}^T(q,p)=(q+h \, p,p),\qquad \phi_{h}^V(q,p)=(q,p-h\nabla V(q)),
\]
but one evaluation of the force per processor and per step can be saved if we consider instead the symmetric composition
\begin{equation} \label{s-v}
   S_h(q,p)=\phi_{h/2}^T\circ\phi_{h}^V\circ\phi_{h/2}^T(q,p),
\end{equation}   
which, for simplicity, we will also refer to as a  St\"ormer--Verlet scheme.

We take $\mu = 1$ and initial conditions
\[
q(0)=(1-e,0),\qquad p(0)=\left(0,\sqrt{(1+e)(1-e)^{-1}}\right), 
\]
so that, if $0 \le e < 1$, the total energy is $H = H_0 = -1/2$, the solution is periodic with period $2 \pi$ and the trajectory is an ellipse of eccentricity $e$.
In the experiments we fix $e = 0.25$ and compute the relative error in phase space and in energy. In the first case, 
given a final time $t_f$ and step size $h = t_f /N$, we compute 
\[
\max_{0\leq n \leq N}\left(\frac{\lVert(q(t_n),p(t_n))-(q_n,p_n)\rVert}{\lVert(q_n,p_n)\rVert}\right), \qquad t_n=n h,
\]
where $(q_n,p_n)$ denotes the numerical approximation, and the exact solution $q(t_n),p(t_n)$ is determined by using the iterative algorithm described in \cite[Sec. 1.5]{blanes16aci} (see also \cite[p. 165]{danby88foc}. In the second case, we compute
\[
 \left\lvert\frac{H_0-H(q_n,p_n)}{H_0}\right\rvert
\]
as a function of time.

\paragraph{Example 2: Lotka-Volterra model.} 

It corresponds to the system of differential equations in $\R^{2}$ given by \cite{hairer06gni}
\begin{equation} \label{lv1}
 u'=u(v-2),\qquad v'=v(1-u),
\end{equation}
and has 
\[
  I(u,v)=\ln u - u + 2 \ln v -v 
\]
as an integral of the motion: $I(u(t), v(t)) = I_0 = \mbox{constant}$ for all $t$. In this case the 
 basic second-order symmetric integrator we use is  $S_h(q,p)=\phi_{h/2}^A\circ\phi_{h}^B\circ\phi_{h/2}^A(q,p)$, with 
\[
\phi_{h}^A(u,v)=(u\exp(h(v-2)),v),\qquad \phi_{h}^B(u,v)=(u,v\exp(h(1-u))).
\]
Here we take as initial conditions $u(0)=1, v(0)=1$, which fix the value of the constant of integration $I(u,v)$ to $I_0=-2$, and integrate up to a final time
$t_f$. Given the number of time steps $N$, we examine the integration error of the numerical solution $(u_n,v_n)$ determined via
\[
\frac{1}{\lfloor0.8N\rfloor}\sum_{n=\lfloor0.8N\rfloor}^N\frac{\lVert (u(t_n),v(t_n))-(u_n,v_n)\rVert }{\lVert(u_n,v_n)\rVert},\qquad t_n=n t_f/N = n h,
\]
where the `exact' solution $u(t_n),v(t_n)$ is determined numerically to high accuracy.
%
On the other hand, the error in $I$ is monitored by evaluating
\[
\left\lvert\frac{I_0-I(u_n,v_n)}{I_0}\right\rvert.
\]

In the graphs collected next, the notation ``\# of evals" is defined as the number of evaluations of the basic scheme $S_h$ per processor (i.e., for $N$ steps of the fourth-order compositions of two-stage integrators, the number of evaluations is $2N$ and so on). In the case of extrapolation methods, it corresponds to the
maximum number of evaluations of $S_h$ required to form the linear combination.

\subsection{Rounding Errors}
\label{sec:Rounding}
Since one sums the outputs of various maps, methods of the form (\ref{eq.3.1o}) can lead to rounding errors which give results for small values of $h$ far worse than could be expected from the theoretical order of the method. One may eliminate or at least reduce such effects in various ways as discussed in \cite{sofroniou08emi}. When designing the method, as mentioned above, one might seek to ensure that the values of $b_i$ are not too large and are of the same scale (e.g. in the range $[0.1,1]$) for example. This offers at least one way to limit the addition and subtraction of numbers which differ greatly in the number of digits after the decimal point (working with a fixed precision). However, a more reliable way to achieve the same goal, at least when the underlying method allows, is to reformulate the algorithms such that in place of working directly with the solutions output by the maps $x^i_h \equiv \prod_{j=1}^mS_{a_{ij}h}(x_0)$ one instead calculates increments for the maps $\Delta x_h^i=x_h^i-x_0$ and then sums the increments $\Delta x_h= \sum_i b_i \, \Delta x_h^i$ before finally calculating $x_h=x_0+\Delta x_h$.

For example, in the case of the St\"ormer--Verlet integrator (\ref{s-v}) used for the Kepler problem, in place of defining $S_h(q,p)$ as the output $x=(q,p$) of the sequence (with $x=(q,p), f(q)=-\nabla V(q)$)
\begin{align}\label{eq:PosVerlet}
q&\gets q+(1/2)hp\nonumber\\
p&\gets p+hf(q)\nonumber\\
q&\gets q+(1/2)hp,
\end{align}
and the composition rule as $(S_h\circ S_h)(x)=S_h(S_h(x))$,
one takes $\widetilde{S}_h(q,p)$ as the output $\Delta x=(\Delta q,\Delta p$) of the sequence
\begin{align}\label{eq:PosVerletInc}
\Delta p&\gets hf(q+(1/2)hp)\nonumber\\
\Delta q&\gets h(p+(1/2)\Delta p)
\end{align}
with the composition rule $\widetilde{S}_h\circ\widetilde{S}_h(x)=\widetilde{S}_h(x+\widetilde{S}_h(x))$. Note that, irrespective of the underlying method, one may always modify the composition rule thusly and define $\widetilde{S}_h(x)=S_h(x)-x$ although with less improvements in accuracy than would be the case otherwise. This has been carried out for the Lotka--Volterra problem, with reasonable success.

With the underlying method and composition redefined to output increments, the sum $\Delta x_h=\sum_ib_i \, \Delta x_h^i$ may be calculated and likewise $x_h=x_0+\Delta x_h$. This form of summation is used in all numerical experiments shown here.

\subsection{Order 4}
\begin{figure}[!htb]
\begin{subfigure}{.5\textwidth}
\includegraphics[width=\textwidth]{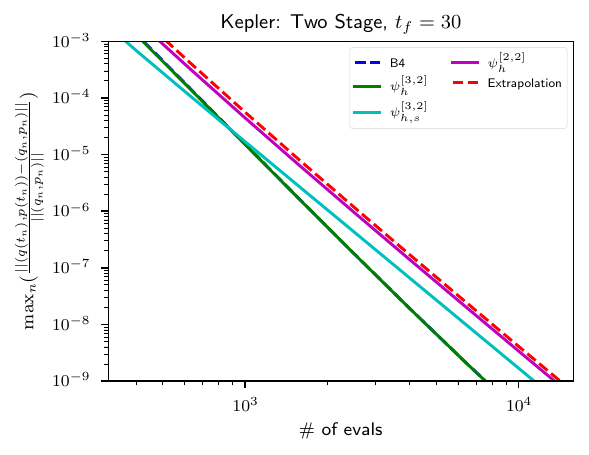}
\caption{}
\end{subfigure}
\begin{subfigure}{.5\textwidth}
\includegraphics[width=\textwidth]{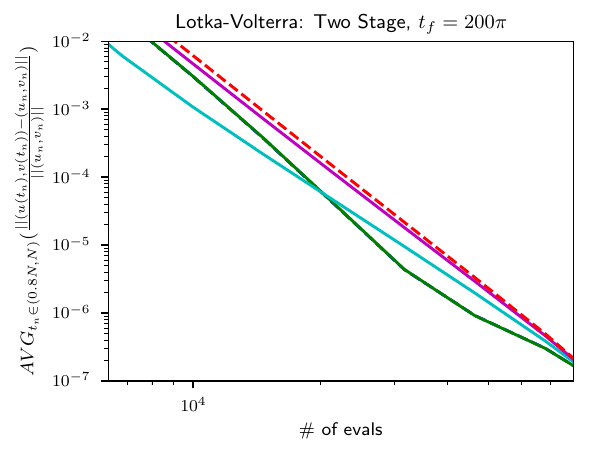}
\caption{}
\end{subfigure}
\begin{subfigure}{.5\textwidth}
\includegraphics[width=\textwidth]{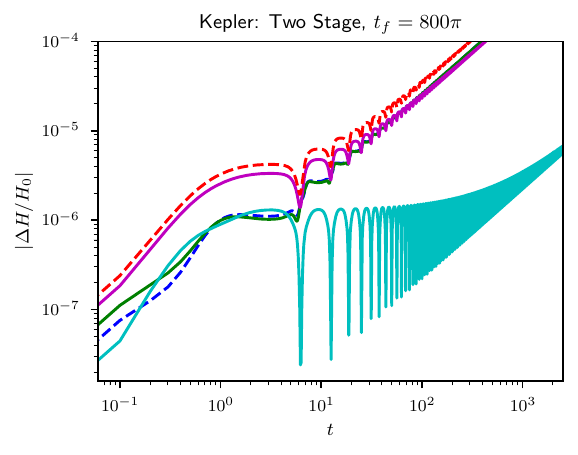}
\caption{}
\end{subfigure}
\begin{subfigure}{.5\textwidth}
\includegraphics[width=\textwidth]{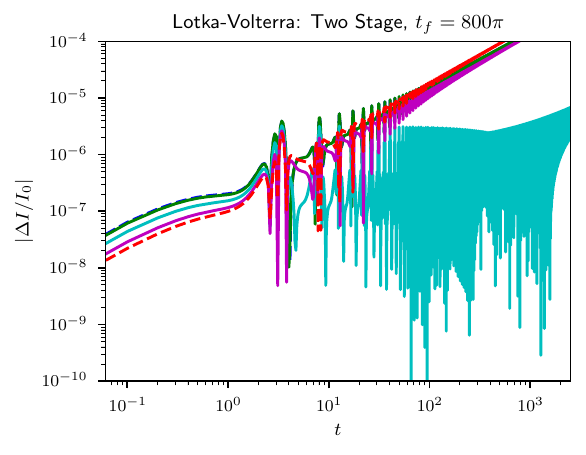}
\caption{}
\end{subfigure}
\caption{\textbf{Order 4}. Efficiency diagrams (a), (b), and relative error in energy vs. time (c), (d) for the Kepler problem and the Lotka--Volterra system. 
All methods have $m=2$, so that they involve 2-stage compositions. In all four plots B4 and $\psi_h^{[3,2]}$, which (numerically) vanish the same terms in the error expansion, are virtually indistinguishable.}
\label{fig:2stage}
\end{figure}

We examine the fourth order methods $\psi_h^{[2,2]}$, $\psi_h^{[3,2]}$ and $\psi_{h,s}^{[3,2]}$ derived in section \ref{sect3}, all of them of type (\ref{psi_or4})
and present comparison results for the two problems, with the fourth-order ($k=2$) extrapolation method and the scheme with $k=3$ derived in \cite{blanes21npi} (called here B4) as the baselines. The corresponding results are collected in Figure \ref{fig:2stage}.

Method $\psi_h^{[2,2]}$ improves slightly on the extrapolation method, in line with expectations, whereas B4 and the optimized scheme $\psi_h^{[3,2]}$ 
(both vanishing the same coefficient $G_{51}$ in the $Z(h)$ expansion) perform equally well in both problems. On the other hand, scheme $\psi_{h,s}^{[3,2]}$,
which is pseudo-symplectic of order 7 (as opposed to $5$ for the other methods shown) performs reasonably well for short time integration - its advantages for long time integration are shown in diagrams (c) and (d) of Figure \ref{fig:2stage}: it begins to show a linear growth in the energy error at times $t$ approximately 100 times larger than the other methods.

\subsection{Order 6}
\begin{figure}[!htb]
\begin{subfigure}{.5\textwidth}
\includegraphics[width=\textwidth]{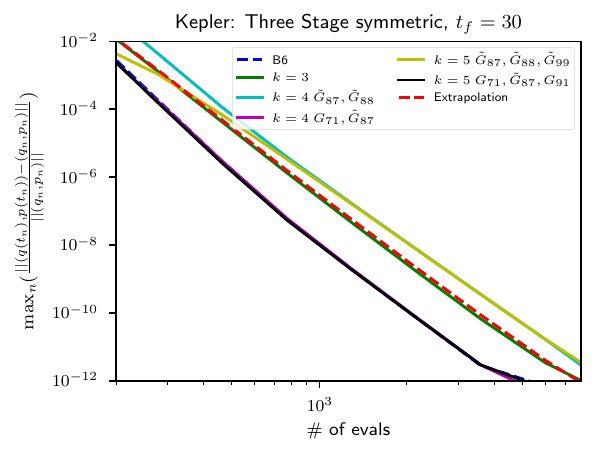}
\caption{}
\end{subfigure}
\begin{subfigure}{.5\textwidth}
\includegraphics[width=\textwidth]{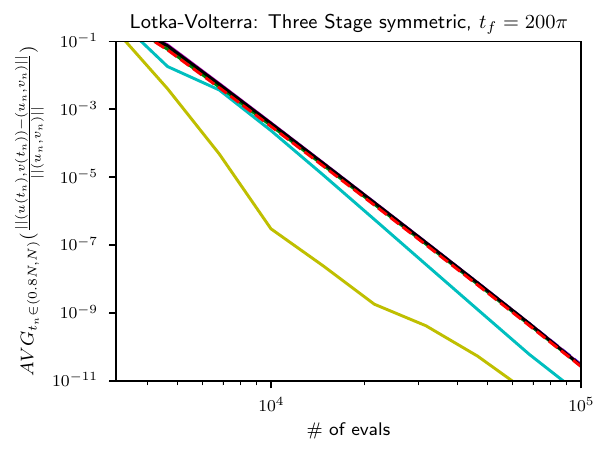}
\caption{}
\end{subfigure}
\begin{subfigure}{.5\textwidth}
\includegraphics[width=\textwidth]{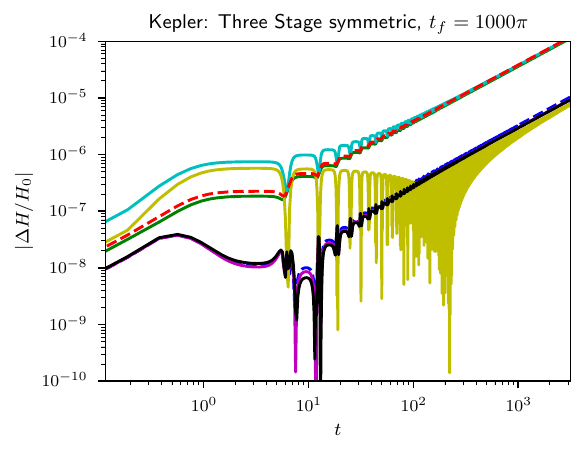}
\caption{}
\end{subfigure}
\begin{subfigure}{.5\textwidth}
\includegraphics[width=\textwidth]{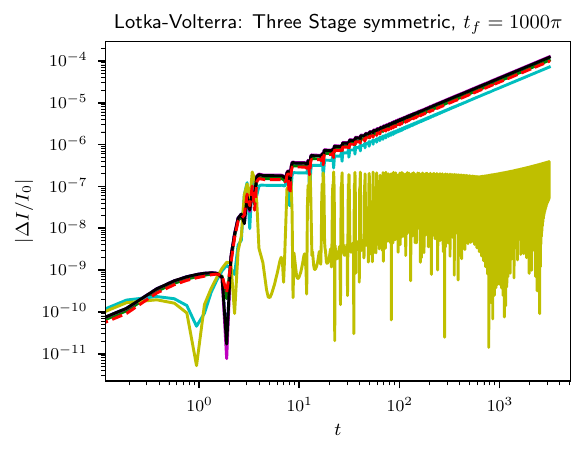}
\caption{}
\end{subfigure}
\caption{\textbf{Order 6}. Efficiency diagrams (a), (b), and relative error in energy vs. time (c), (d) for the Kepler problem and the Lotka--Volterra system. 
All methods have $m=3$, and they involve 3-stage symmetric compositions. In all four plots B6 and the $k=5$ integrator, which both vanish $G_{71},\tilde{G}_{87},G_{91}$, 
are virtually indistinguishable.}
\label{fig:3stagesymm}
\end{figure}

We first examine the sixth order methods derived in section \ref{sect3} and involving 3-stage symmetric compositions, i.e, linear combinations of the form
(\ref{psi_or6}), and present comparison results for the two problems, with the sixth-order ($k=3$) extrapolation method and the $k=5$ scheme proposed in
\cite{blanes21npi} (called here B6) as the baselines. Figure \ref{fig:3stagesymm} collects the relative errors in position and the integrals of motion for both
problems. 

As for the two stage case, the method with $k=3$ we have constructed improves only very slightly on the extrapolation method $\psi_h^{[6]}$ of eq. (\ref{em46}).
The schemes which eliminate $G_{71}$ (which includes method B6) are clearly superior in the Kepler problem, but perform the same in the Lotka--Volterra
system. Methods of higher pseudo-symplectic order perform worse in the first case, but show massive improvements in the second: in fact, one is able to achieve parity with B6 using only 4 processors. Our proposed $k=5$ scheme which vanishes terms $\tilde{G}_{87}$, $\tilde{G}_{88}$, and $\tilde{G}_{99}$, and thus is 
pseudo-symplectic of order 9, shows the best performance for the Lotka--Volterra problem, and in addition provides improved conservation of the first integrals in
both problems.

\begin{figure}[h]
\begin{subfigure}{.5\textwidth}
\includegraphics[width=\textwidth]{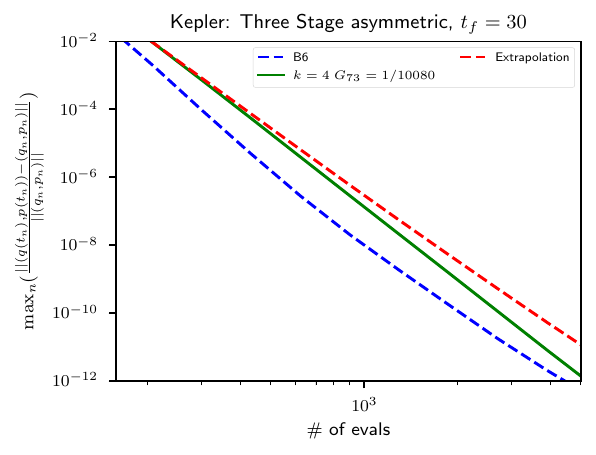}
\caption{}
\end{subfigure}
\begin{subfigure}{.5\textwidth}
\includegraphics[width=\textwidth]{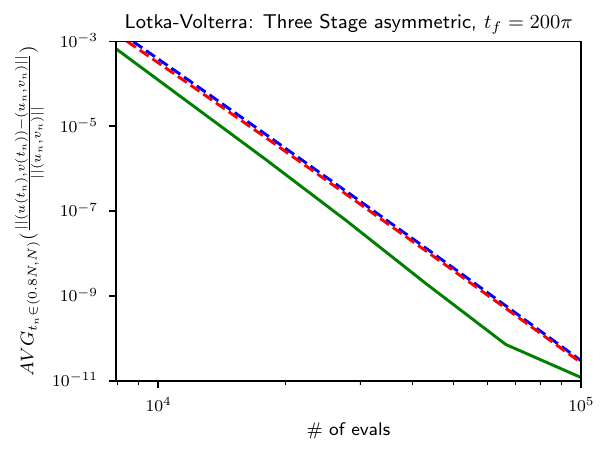}
\caption{}
\end{subfigure}
\caption{\textbf{Order 6}. Efficiency diagrams for both problems. The proposed scheme based on a linear combination of $k=4$ non-palindromic compositions
shows an improved behavior on the Lotka--Volterra system.}
\label{fig:6asymmetric}
\end{figure}

On the other hand, the method constructed in section \ref{sect3} 
involving non-palindromic compositions \eqref{npc6} with $k=4$,
whilst not improving on B6 in the Kepler problem, owing to the $\tilde{G}_{87}$ term, does show improved performance in the Lotka-Volterra problem,
as shown in Figure \ref{fig:6asymmetric}.

\subsection{Order 8}
\begin{figure}[h]
\begin{subfigure}{.5\textwidth}
\includegraphics[width=\textwidth]{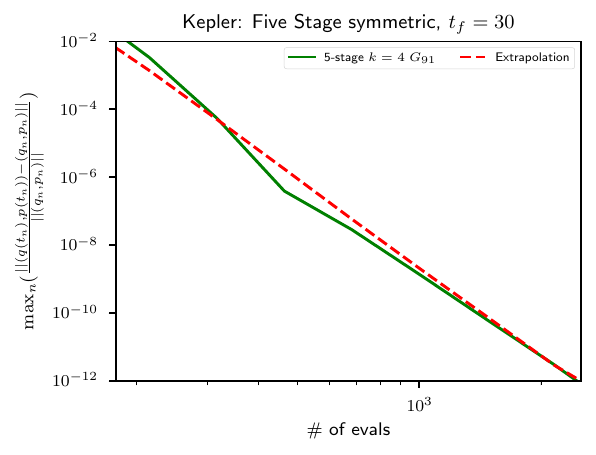}
\caption{}
\end{subfigure}
\begin{subfigure}{.5\textwidth}
\includegraphics[width=\textwidth]{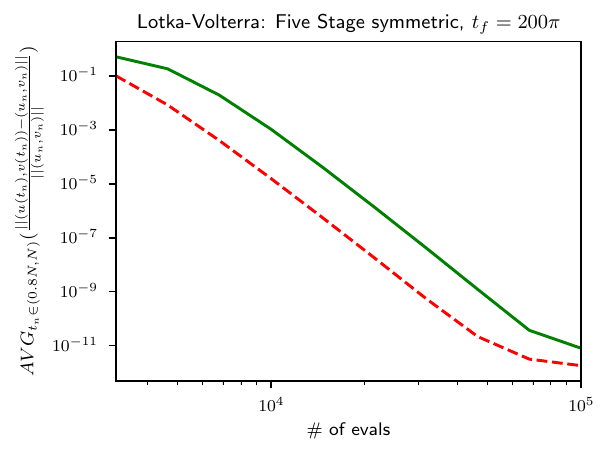}
\caption{}
\end{subfigure}
\caption{\textbf{Order 8}. Efficiency diagrams for both problems. Including an additional stage in the composition renders a better performance with respect to the standard
MPE only for the Kepler problem.}
\label{fig:5stage}
\end{figure}

We compare the integrator proposed in section \ref{sect3} (a linear combination of $k=4$ symmetric compositions involving 5 stages) with the extrapolation method $\psi_h^{[8]}$, which uses at most 4 evaluations of $S_h$ per processor. Now, the elimination of the $G_{91}$ term by the new scheme 
only compensates the extra cost of the 
additional evaluation of $S_h$ in the Kepler problem (see Figure \ref{fig:5stage}). This is so in part because the extrapolation method $\psi_h^{[8]}$ works better
than expected only from local error considerations: in fact, the order exhibited in both diagrams of Figure \ref{fig:5stage} is $\approx 8.4$.

\section{Low-latency methods}
\label{sect5}

As stated in the introduction, \emph{latency} may considerably affect the overall performance of the standard MPE (\ref{eq.extrap}) when it is implemented
in a parallel environment, so that it is relevant to analyze the situation where the summation of the different compositions 
is delayed by, say, $p$ steps. In other words, instead of (\ref{linear_comb1}) one has
linear combinations of the form (\ref{linear_comb_p}):
\begin{equation*}
\psi^{[k,m,p]}_h=\sum_{i=1}^k b_i\left(\prod_{j=1}^mS_{a_{ij}h}\right)^p, \qquad\quad \sum_{i=1}^kb_i=\sum_{j=1}^ma_{ij}=1,
\end{equation*}
which is naturally a special case of (\ref{linear_comb1}) with $\bar{m}=mp$ and a periodicity restriction on the sequence
\begin{equation}\label{eq:LowLatencyNovelInteg}
\psi^{[k,mp]}_h=\sum_{i=1}^k b_i\prod_{j=1}^{mp}S_{a_{ij}h}, \qquad\quad \sum_{i=1}^kb_i=\sum_{j=1}^{mp}a_{ij}=1,\qquad a_{i,j}=a_{i,j+m}.
\end{equation}
Relaxing the restrictions on $a_{ij}$ in (\ref{eq:LowLatencyNovelInteg}) could generate very high-order methods:  for illustration, 
even with a symmetry restriction, for $mp=2\times 10=20$ and $k=3$ processors one has $29$ free parameters, which in principle would enable one to 
satisfy the 22 order conditions required by a 9th-order method, with $7$ free parameters for optimization. However, solving 
such a system is a difficult computational challenge, and of course then enforces multiples of $20$ evaluations if one desires to further delay summation.

Since this delayed sum presents obvious computational advantages in a parallel environment, it is indeed relevant to analyze how the error of the
corresponding method
depends on $p$, and in particular, whether there is some order reduction
when one uses very large $p=\mathcal{O}(h^{-1})$. If we denote by $\Psi^{[k,m,p]}(h)$ the operator associated to the linear combination $\psi^{[k,m,p]}_h$, then
a similar calculation as was done in section \ref{sect3} leads us to
\[
  \Psi^{[k,m,p]}(h) = \exp(p h Y_1/2) \, Z_p(h) \, \exp(p h Y_1/2),
\]
with  
\begin{equation} \label{expan_Zp}
\begin{aligned}
  & Z_p(h) = G_{00} \, I + p h^3  G_{31} E_{31} + p h^4 G_{41} E_{41}  + p h^5 
   \left( \sum_{\ell =1}^2 G_{5 \ell} \, E_{5 \ell} + \frac{1}{24} (p^2-1) G_{31} E_{52} \right)\\
  & \qquad + p h^6 \left( \sum_{\ell =1}^2 G_{6 \ell} E_{6 \ell} + \frac{1}{24} (p^2-1) G_{41} E_{62}  + \frac{1}{2} \,p \, \tilde{G}_{63} E_{31}^2 \right) + \mathcal{O}(ph^7)
\end{aligned}   
\end{equation}
If the method is of order 4, then $G_{00} = 1$, $G_{31}=G_{41}=0$, so that, after $n$ steps, one has
\begin{equation} \label{aftern}
  \left( \Psi^{[k,m,p]}(h) \right)^n = \exp(V_n(h)) + \mathcal{O}(h^{10})
\end{equation}
with
\[
\begin{aligned}
  & V_n(h) = n p h Y_1 + n p h^5 \sum_{\ell =1}^2 G_{5 \ell} \, E_{5 \ell} + n p h^6 \left( \sum_{\ell =1}^2 G_{6 \ell} E_{6 \ell}  
  + \frac{1}{2} \,p \, \tilde{G}_{63} E_{31}^2 \right)  \\
  & \qquad\quad +  n p h^7 \left( \sum_{\ell =1}^4 G_{7 \ell} E_{7 \ell} +   \frac{1}{2} \,p \, \tilde{G}_{75} \{E_{31}, E_{41} \} \right)
 + \mathcal{O}(n ph^8).
\end{aligned} 
\]
To arrive at the final time $t_f$ after $N$ steps, $t_f = N h$, we need to take $n = N/p$ in (\ref{aftern}), so that
\[
 V_{N/p}(h) = t_f Y_1 + t_f h^4 \sum_{\ell =1}^2 G_{5 \ell} \, E_{5 \ell} + t_f h^5 \sum_{\ell =1}^2 G_{6 \ell} E_{6 \ell} + 
 t_f \, p \, h^5   \frac{1}{2} \tilde{G}_{63} E_{31}^2 + 
 \mathcal{O}(t_f h^6).
\]
The error thus grows with $p$, and in the extreme case where one performs the sum once so that $p=t_f/h$, the method has global error $\mathcal{O}(h^4+t_f h^4)$.
If, in addition, the error terms $G_{5j}, E_{5j}$ vanish or are relatively small, then we have an apparent order reduction. This is the case of methods B4 and
$\psi_{h}^{[3,2]}$ for the Kepler problem, since $G_{51}=0$ and $G_{52}, E_{52}$ are very small (cf. Figure \ref{fig:2stage}).

\begin{figure}[h]
\begin{subfigure}{.5\textwidth}
\includegraphics[width=\textwidth]{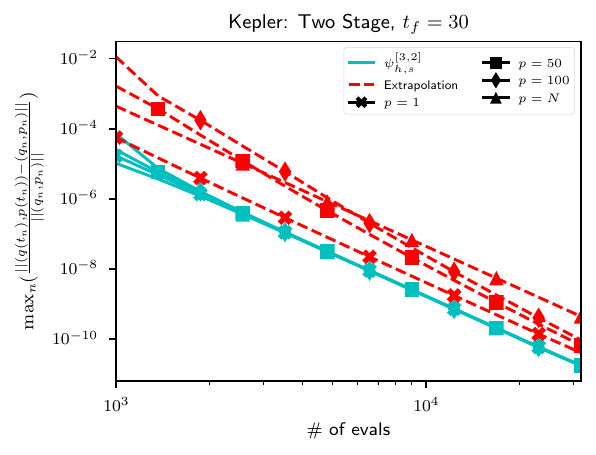}
\caption{}
\end{subfigure}
\begin{subfigure}{.5\textwidth}
\includegraphics[width=\textwidth]{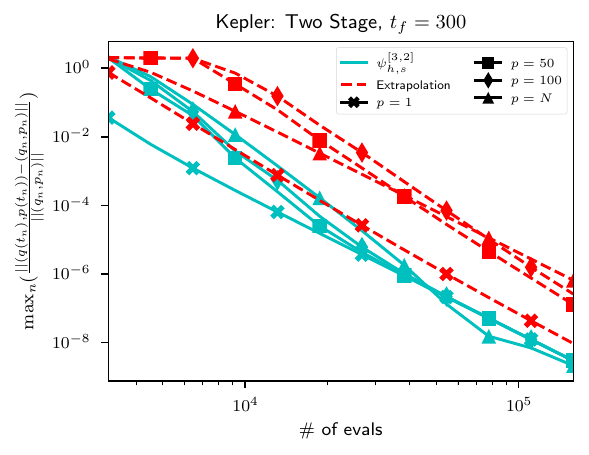}
\caption{}
\end{subfigure}
\caption{\textbf{Order 4} when the sum is computed after $p$ steps. Extrapolation method and scheme $\psi_{h,s}^{[3,2]}$. The results obtained
by the pseudo-symplectic method of order 7 do not depend on the value of $p$, in accordance with the analysis. The advantages 
become more apparent on increasing the final integration time $t_f$ (right panel).}
\label{fig:TwostageLatency}
\end{figure}

On the other hand, for the 4th-order scheme $\psi_{h,s}^{[3,2]}$ (pseudo-symplectic of order 7), it holds that both $\tilde{G}_{63}$ and $\tilde{G}_{75}$ vanish,
so that the global error is constant with $p$ and $\mathcal{O}(h^4)$ regardless of $p$ and $t_f$. This method is compared to the (representative) fourth order extrapolation method $\psi_h^{[4]}$ in Figure \ref{fig:TwostageLatency}. We clearly see that $\psi_{h,s}^{[3,2]}$ performs the same regardless of the value of $p$, even if the sum is done only once, at the end of the integration $p=N$. In contrast, other methods will suffer decreasing performance and even (effective) order reduction for large $p$. The benefits of preserving the symplectic character up to a higher order become more apparent on increasing the final integration time $t_f$, see Figure \ref{fig:TwostageLatency} (b).

\begin{figure}[h]
\centering
\includegraphics[width=0.8\textwidth]{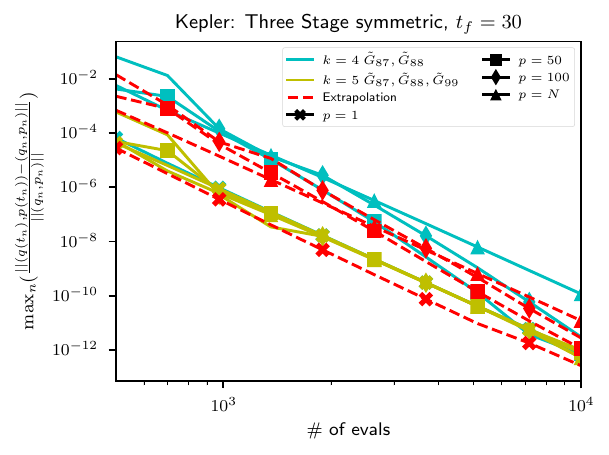}
\caption{\textbf{Order 6}. Comparison of the performance of different schemes with $p$. Only the method which vanish $\tilde{G}_{87}$, $\tilde{G}_{88}$ and $\tilde{G}_{99}$ (pseudo-symplectic of order 9) exhibits an $\mathcal{O}(h^6)$ error, independently of $p$.}
\label{fig:ThreestagesymmLatency}
\end{figure}

In the case of a sixth-order method built with symmetric compositions, the corresponding operator $V_n(h)$ in (\ref{aftern}) reads
\[
\begin{aligned}
  & V_n(h) =  n p h Y_1 + n p h^7 \, \sum_{\ell =1}^4 G_{7 \ell} E_{7 \ell}   
+  n p^2 h^8 \frac{1}{2} \left(  \tilde{G}_{87} \{ E_{31}, E_{51} \}  +
      \tilde{G}_{88} \{ E_{31}, E_{52} \} \right)  \\
      & \qquad\quad + n p^3 h^9 \frac{1}{6} \tilde{G}_{99} + \mathcal{O}(n p h^9)
\end{aligned}      
\]
so that at the final time $t_f$ we have
\[
 V_{N/p}(h) = t_f Y_1 + t_f h^6\sum_{\ell=1}^4G_{7\ell}E_{7\ell}+\frac{t_f ph^7}{2}\left(\tilde{G}_{87}\tilde{E}_{87}+\tilde{G}_{88}\tilde{E}_{88}\right)+\frac{t_f p^2h^8}{6}\tilde{G}_{99}\tilde{E}_{99}+\mathcal{O}(t_f h^7+t_f p^2h^9). 
\]
Notice that both the extrapolation $\psi_h^{[6]}$ and B6 methods exhibit a $\mathcal{O}(h^6+ph^7+p^2h^8)$ global error, which in the worst case 
($p=t_f/h$) gives $\mathcal{O}(t_f^2h^6)$, whereas the $k=4$ method constructed in section \ref{sect3} that vanishes $\tilde{G}_{87}$, $\tilde{G}_{88}$
possesses the same $\mathcal{O}(t_f^2h^6)$ dependence in the global error. Finally, the $k=5$ scheme vanishing the terms $\tilde{G}_{87}$, $\tilde{G}_{88}$ and $\tilde{G}_{99}$ has an $\mathcal{O}(h^6)$ error. All these features are clearly visible in Figure \ref{fig:ThreestagesymmLatency}: 
the performance of methods which fail to vanish $\tilde{G}_{99}$ always degrades upon increasing $p$ (although for sufficiently small $h$, i.e. large number of 
evaluations, one is able to make the higher order non-symplectic terms small enough to improve performance). Notice that 
the method of highest pseudo-symplectic order shows invariable performance, as was the case for the 4th-order method in Figure \ref{fig:TwostageLatency}.
We see then that, by choosing parameters in \eqref{linear_comb1} so as to increase the order of structure-preservation of the method, it is possible to delay the sum by
a number of steps (and therefore reducing the latency problem) without affecting the overall performance of the scheme.

\subsection*{Acknowledgments}
The authors are grateful to A. Escorihuela-Tom\`as for solving the equations leading to the methods of order 8 reported in this paper. 
This work has been
funded by Ministerio de Ciencia e Innovaci\'on (Spain) through project PID2022-136585NB-C21, 
MCIN/AEI/10.13039/501100011033/FEDER, UE, and also by Generalitat Valenciana (Spain) through project CIAICO/2021/180.

\subsection*{Competing Interests}
The authors declare no competing interests.

\section*{Appendix: Coefficients}
We collect here the coefficients of three of the most efficient methods we have constructed.

\begin{table}[h!]
\centering
 \caption{Coefficients for the method $\psi_{h,s}^{[3,2]}$ in Figure \ref{fig:2stage}, of the type shown in Equation \eqref{psi_or4}.}
 \begin{tabular}{||l|l||} 
 \hline
 $a_1=-0.19220568886474299$   &$b_1=0.09012936855999465 $\\
 $a_2=0.7952090547057717$  & $b_2=-1.8742613286568583$\\
 $a_3=0.615$ & $b_3=1-b_1-b_2$\\
 \hline
 \end{tabular}
\end{table}

\begin{table}[h!]
\centering
\caption{Coefficients for the 6th-order method with $k=5$ which
vanishes the terms $\tilde{G}_{87}$, $\tilde{G}_{88}$, and $\tilde{G}_{99}$ in Figure \ref{fig:3stagesymm}, of the type shown in Equation \eqref{psi_or6}.}
 \begin{tabular}{||l|l||} 
 \hline
 $a_1=0.7702669932516844$   &$b_1= 0.7482993205697204 $\\
 $a_2=2/100$  & $b_2=-0.34096002148336635$\\
 $a_3=0.5133170199053506$ & $b_3= -1.5697387622875072$\\
 $a_4= 1.1686905913031624 $  & $b_4= -0.11572553679884676$\\
 $a_5=1/3$  & $b_5=1-b_1-b_2-b_3-b_4$\\
 \hline
 \end{tabular}
\end{table}

\begin{table}[h!]
\centering
\caption{Coefficients for the 8th-order method with $k=4$ which
vanishes $G_{91}$ in Figure \ref{fig:5stage}, of the type shown in Equation \eqref{psi_or8}.}
 \begin{tabular}{||l|l|l||} 
 \hline
 $a_{11}=-0.2539842055534987$ &  $a_{12}= 0.4514159659747628$ & $b_1= 0.6402721677360648 $\\
 $a_{21}=-0.1297472147351918$ & $a_{22}= 0.5893868250930246$ & $b_2= -0.4488395035838362$\\
 $a_{31}=0.283267969084071$ & $a_{32}=  0.0411275969512266$ & $b_3= -11.611098146500447$\\
 $a_{41}= 0.0671551220219572$ &  $a_{42}=  0.3228966120312048$& $b_4=1-b_1-b_2-b_3$\\
 \hline
 \end{tabular}
\end{table}

\bibliographystyle{siam}

\end{document}